\documentclass[11pt]{amsart}
\usepackage{amscd,amssymb}
\usepackage[arrow,matrix,frame,poly,arc]{xy}

\theoremstyle{plain}
\newtheorem{theorem}[subsection]{Theorem}
\newtheorem{lemma}[subsection]{Lemma}
\newtheorem{prop}[subsection]{Proposition}
\newtheorem{corollary}[subsection]{Corollary}

\newtheorem{thm}{Theorem}

\theoremstyle{definition}
\newtheorem{remark}[subsection]{Remark}
\newtheorem{definition}[subsection]{Definition}
\newtheorem{example}[subsection]{Example}

\newtheorem*{ack}{Acknowledgments}

\numberwithin{equation}{section}

\newcommand{\A}{{\mathcal A}}
\newcommand{\G}{{\mathcal G}}
\newcommand{\I}{{\mathcal I}}
\newcommand{\J}{{\mathcal J}}

\newcommand{\LL}{{\mathbf L}}
\newcommand{\FF}{{\mathbf F}}

\newcommand{\Z}{\mathbb{Z}}
\newcommand{\Q}{\mathbb{Q}}

\newcommand{\C}{\mathbb{C}}
\newcommand{\F}{\mathbb{F}}
\newcommand{\K}{\Bbbk}

\newcommand{\B}{{\mathfrak B}}
\newcommand{\Amod}{{\mathfrak A}}
\newcommand{\HH}{{\mathfrak H}}
\newcommand{\m}{{\mathfrak m}}
\renewcommand{\r}{{\overline{r}}}

\DeclareMathOperator{\Hom}{Hom}
\DeclareMathOperator{\Hilb}{Hilb}
\DeclareMathOperator{\rank}{rank}
\DeclareMathOperator{\gr}{gr}
\DeclareMathOperator{\im}{im}
\DeclareMathOperator{\coker}{coker}

\DeclareMathOperator{\id}{id}

\DeclareMathOperator{\pr}{pr}
\DeclareMathOperator{\lk}{lk}

\DeclareMathOperator{\Sym}{Sym}

\DeclareMathOperator{\ch}{char}

\newcommand{\surj}{\twoheadrightarrow}

\newcommand{\abs}[1]{\left| \, #1 \right|}

\newenvironment{alphenum}{
\begin{enumerate}}{\end{enumerate}}

\begin{document}

\title[Chen Lie algebras]{Chen Lie algebras}

\author[Stefan Papadima]{Stefan Papadima}
\address{Institute of Mathematics of the Academy,
P.O. Box 1-764,
RO-70700 Bucharest, Romania}
\email{Stefan.Papadima@imar.ro}

\author[Alexander~I.~Suciu]{Alexander~I.~Suciu}
\address{Department of Mathematics,
Northeastern University,
Boston, MA 02115, USA}
\email{a.suciu@neu.edu}
\urladdr{http://www.math.neu.edu/\~{}suciu}

%



\begin{abstract}
The Chen groups of a finitely-presented group $G$ are the lower central
series quotients of its maximal metabelian quotient, $G/G''$. The direct
sum of the Chen groups is a graded Lie algebra, with bracket induced by
the group commutator. If $G$ is the fundamental group of a formal space,
we give an analog of a basic result of D.~Sullivan, by showing that the rational 
Chen Lie algebra of $G$ is isomorphic to the rational holonomy Lie 
algebra of $G$ modulo the second derived subalgebra. Following an 
idea of W.S.~Massey, we point out a connection between the Alexander 
invariant of a group $G$ defined by commutator-relators, and its integral 
holonomy Lie algebra.

As an application, we determine the Chen Lie algebras of several classes of 
geometrically defined groups, including surface-like groups, fundamental 
groups of certain link complements in $S^3$, 
and fundamental groups of complements 
of hyperplane arrangements in $\C^{\ell}$.  For link groups, we sharpen 
Massey and Traldi's solution of the Murasugi conjecture. For arrangement 
groups, we prove that the rational Chen Lie algebra is combinatorially determined.
\end{abstract}
\maketitle

\section{Introduction}
\label{sect:intro}

\subsection{}
A classical construction of W.~Magnus associates to a group $G$ 
a graded Lie algebra over $\Z$, 
\[
\gr(G)=\bigoplus\nolimits_{k\ge 1} \Gamma_k G/\Gamma_{k+1} G,
\]
where $\{\Gamma_k G\}_{k\ge 1}$ is the lower central series  
of the group, defined inductively by $\Gamma_1 G=G$ and 
$\Gamma_{k+1} G=(\Gamma_{k} G,G)$, and the Lie bracket $[x,y]$ 
is induced from the group commutator $(x,y)=xyx^{-1}y^{-1}$.  

Many properties of a group are reflected in properties of its associated 
graded Lie algebra.  For instance, if $G$ is finitely generated, 
then the abelian groups $\gr_k(G)$ are also finitely generated;  
their ranks, $\phi_k(G)$, are important numerical invariants of $G$. 
In the case when $G=\FF_n$, the free group of rank $n$, 
Magnus showed that $\gr(\FF_n)=\LL_n$, the free Lie algebra 
on $n$ generators, whose graded ranks were computed by E.~Witt. 
In general, though, the computation of the LCS ranks $\phi_k(G)$ 
can be exceedingly difficult.  

\subsection{}
In his thesis \cite{Ch51}, K.T.~Chen introduced a more manageable 
approximation to the LCS ranks. Let $G/G''$ be the maximal metabelian 
quotient of $G$. The Chen groups are the graded pieces of 
the associated graded Lie algebra $\gr(G/G'')$. 
Assume $G$ is finitely-generated, and let $\theta_k (G)=\rank (\gr_k(G/G''))$ 
be the rank of the $k$-th Chen group. 
Then $\theta_k(G)=\phi_k(G)$, for $k\le 3$, and
$\theta_k(G)\le \phi_k(G)$, for $k> 3$. 

K.T.~Chen showed that the Chen groups of $\FF_n$ are torsion-free, and 
computed their ranks. He also gave an algorithm for computing the ranks 
$\theta_k(G)$ for an arbitrary finitely presented group $G$, but that algorithm is  
highly impractical. 

\subsection{}
In a subsequent paper \cite{Ch}, Chen introduced the
{\em rational holonomy Lie algebra} of a space $X$. 
Assuming $X$ is a connected CW-complex, with finite $2$-skeleton, 
this Lie algebra is 
\[
\HH (X;\Q) := \LL^*(H_1(X; \Q))/\, \text{ideal} \big( \im 
(\partial^{\Q}_X)\big),
\]
where $ \LL^*(H_1(X; \Q))$ is the free (graded) Lie algebra over $\Q$, 
generated in degree $1$ by $H_1(X; \Q)$, and $\partial^{\Q}_X$ is
the dual of the cup-product map,
$\cup^{\Q}_X\colon H^1(X; \Q) \wedge H^1(X; \Q) \to H^2(X; \Q)$.  
It is readily seen that $\HH(X;\Q) =\HH(G;\Q)$, where $G$ is the 
fundamental group of $X$, and $\HH(G;\Q)$ is the 
rational holonomy Lie algebra of the 
Eilenberg-MacLane space $K(G,1)$. 

\subsection{}
Now suppose $X$ is a {\em formal} space, in the sense of Sullivan~\cite{Su}.
Loosely speaking, this means that the rational homotopy type of $X$ is a
formal consequence of the rational cohomology algebra of $X$.
Let  $G=\pi_1(X)$.  Sullivan showed that
there is an isomorphism of graded Lie algebras over $\Q$, 
\begin{equation}
\label{sulli=eq}
\gr(G) \otimes \Q \cong \HH (X; \Q) .
\end{equation}

There are many examples of formal spaces, for instance, spheres, 
compact K\"{a}hler manifolds (cf.~\cite{DGMS}), 
and complements of complex hyperplane arrangements (as follows from~\cite{Br}). 
Formality is preserved by taking products or wedges. 

As proved in \cite{Su} and \cite{Mo}, the fundamental group of
a formal space is $1$-formal, in the sense of Definition \ref{def:1formal}. 
The isomorphism $\gr(G) \otimes \Q \cong \HH (G; \Q)$ holds for 
arbitrary $1$-formal groups. 
Fundamental groups of complements of complex projective hypersurfaces 
are always $1$-formal (this was proved by Kohno \cite{Ko83}, using 
resolution of singularities to reduce to the case of complements of
normal-crossing divisors). In general though, complements of normal-crossing
divisors may not be formal, see Morgan \cite{Mo}.

\subsection{}
Assuming $H_1(G; \Z)$ is torsion-free, one may define an integral 
form of Chen's holonomy Lie algebra.  This graded $\Z$-Lie algebra, 
denoted by $\HH(G)$, comes equipped with a natural  
epimorphism of graded Lie algebras, $\Psi_G \colon \HH(G) \surj \gr (G)$.  
If $G$ is $1$-formal and $\HH(G)$ is torsion-free, 
the map  $\Psi_G$ gives an isomorphism
\begin{equation}
\label{jem=eq}
 \HH(G) \cong \gr (G), 
\end{equation}
see~\cite{MP}, and also Section~\ref{zhol=sec}.

\subsection{}
The  purpose of this note is to produce analogs of isomorphisms 
\eqref{sulli=eq} and \eqref{jem=eq} for the Chen Lie algebra $\gr (G/G'')$.  
Without much more effort, we will do this for all the higher-order 
Chen Lie algebras, 
\[
\gr\big(G/G^{(i)}\big), 
\]
$i\ge 2$, where $\{G^{(i)}\}_{i\ge 0}$ is the derived series of $G$, 
defined inductively by $G^{(0)}=G$ and $G^{(i+1)}=(G^{(i)},G^{(i)})$.  

In Sections \ref{sect:malcev}--\ref{sect:holoformal}, we treat the rational case. 
Our first main result is Theorem~\ref{thm:Malcompl}, where we describe
the Malcev completion of $G/G^{(i)}$ by means of a functorial formula, 
in terms of the Malcev completion of $G$.  Under a formality assumption, 
we deduce in Theorem~\ref{thm:hatformal} the following analog of Sullivan's isomorphism~\eqref{sulli=eq}.

\begin{thm}
\label{thm1}
Let $G$ be a finitely presented group. 
 If $G$ is $1$-formal, then there is an isomorphism of graded Lie algebras, 
\begin{equation}
\label{chsu=eq}
\gr \!\big( G/G^{(i)} \big)\otimes \Q \cong \HH (G; \Q)/\HH^{(i)} (G; \Q), \quad 
\forall i \ge 2.
\end{equation}
\end{thm}
Formula~\eqref{chsu=eq} says that the rational associated graded Lie
algebra of the $i$-th derived quotient of the group is isomorphic to the 
$i$-th derived quotient of the $\Q$-holonomy Lie algebra of the group. 
This formula does not hold for arbitrary $G$, see Example~\ref{notf=ex}.

\subsection{}
In Section~\ref{zhol=sec}, we treat the integral case.  
We first show that the map $\Psi_G \colon \HH(G) \surj \gr (G)$ induces 
a natural epimorphism of graded Lie algebras,
$\Psi_G^{(i)} \colon \HH(G)/\HH(G)^{(i)} \surj \gr (G/G^{(i)})$, for each $i\ge 1$.
Under some torsion-freeness and formality hypotheses, 
we deduce in Proposition~\ref{tors} the following analog of isomorphism 
\eqref{jem=eq}.

\begin{thm}
\label{thm2}
Let $G$ be a finitely presented group, with torsion-free abelianization. 
Suppose  $G$ is $1$-formal.  If $\HH(G)/\HH(G)^{(i)}$ is torsion-free,  
then the map  $\Psi_G^{(i)}$ gives an isomorphism 
of graded Lie algebras over $\Z$,
\begin{equation}
\label{jemi=eq}
\HH(G)/\HH(G)^{(i)} \cong \gr (G/G^{(i)}).  
\end{equation}
\end{thm}

\subsection{}
Our approach works particularly well for  commutator-relators groups,
that is, for finitely presented groups $G$ defined by relators belonging
to $\Gamma_2$.  In this case, the integral holonomy Lie algebra $\HH(G)$ 
may be described directly in terms of the defining relations of $G$, see 
Proposition~\ref{mod3=prop}. 

The adjoint representation of $\HH(G)$ endows the 
{\em infinitesimal Alexander invariant}, $B(\HH(G)) = \HH(G)'/\HH(G)''$, 
with a natural graded module structure over the polynomial
ring $S=\Sym(H_1G)$. In Theorem~\ref{thm:ppres}, we derive from 
the defining presentation of $G$ a finite presentation of 
$B(\HH(G))$ over the ring $S$, while in Proposition~\ref{prop:linalex}, 
we relate the module $B(\HH(G))$ to the classical Alexander invariant 
of the group, $B_G=G'/G''$.

\subsection{}
In Theorem~\ref{1rel=thm}, we determine the Chen Lie algebras of 
``surface-like" groups.  If  $G=\langle x_1, \dots, x_n \mid r\rangle$
is $1$-formal, with relator $r\in \Gamma_2 \FF_n$ such that 
$\r = r \bmod \Gamma_3 \FF_n$  
satisfies certain non-degeneracy conditions, then:
\[
\gr (G/G'') = \LL_n/ \big( {\rm ideal} (\r) +\LL''_n\big). 
\]
Moreover, $\gr (G/G'')$ is torsion-free, with Hilbert  series 
$1+ nt - (1-nt+t^2)/(1-t)^n$. 
This result is a Chen Lie algebra analog of work by Labute~\cite{L70}, 
where the Lie algebra $\gr(G)$ was determined for certain 
one-relator groups $G$. 

There is a rich variety of one commutator-relator groups $G$ satisfying the
hypotheses of our Theorem~\ref{1rel=thm}. Among them, we mention the class of
fundamental groups of irreducible complex plane projective curves of
positive genus (see Example~\ref{curve=ex}), and the algebraically defined
class of surface-like groups from Example~\ref{psurf=ex}.

\subsection{}
Let $K=(K_1,\dots ,K_n)$ be a link in $S^3$, with complement $X$,
and fundamental group $G=\pi_1(X)$.
The Chen groups $\gr_k(G/G'')$, first considered by Chen
in  \cite{Ch51}, \cite{Ch52}, were intensively studied in
the 70's and 80's, see  for example \cite{L}, \cite{Md},
\cite{Ms}, \cite{MT}, \cite{Mu},  \cite{Tr89}.
The most complete computation of Chen groups of links was
done by Murasugi \cite{Mu}, in the case of $2$-component links.
No general formula is known for the Chen groups of an arbitrary link.
Nevertheless, Murasugi proposed a simple
formula for the Chen ranks $\theta_k$ and LCS ranks $\phi_k$ of a link group,
in the case when all the linking numbers are equal to $\pm 1$.
Murasugi's conjecture was proved in particular cases by
Kojima \cite{Kj} and Maeda \cite{Md},
and in full generality by Massey and Traldi \cite{MT}. Further work in
this direction  was done in \cite{L}, \cite{MP}, and \cite{BP1}.

As an application of our results, we prove an even stronger
form of Murasugi's conjecture, both rationally and integrally.
The main contribution is at the multiplicative level of the
Chen Lie algebras of the links under consideration. See
Theorem~\ref{thm:qmura} (f), and Theorem~\ref{thm:zmura}
\eqref{zmb}--\eqref{zmc} respectively.

\subsection{}
Let $\A=\{H_1,\dots ,H_n\}$ be an arrangement of hyperplanes in
$\C^{\ell}$, with complement $X$, and group $G=\pi_1(X)$.
The study of lower central series quotients of arrangement groups
was initiated by Kohno~\cite{Ko85}, who computed the LCS ranks of the 
pure braid groups.  Falk and Randell \cite{FR} extended Kohno's 
computation to the broader class of fiber-type arrangement groups, 
expressing the LCS ranks in terms of the exponents of the arrangement.
See \cite{JP}, \cite{PY} for further generalizations of the LCS formula, 
and \cite{S}, \cite{SS} for other formulas, and computations in low ranks, 
even when the LCS formula does not hold.

Another direction was started in \cite{CS1}, \cite{CS2}, with the study 
of Chen groups of arrangements.  The Chen ranks $\theta_k(G)$ 
can provide stronger information than the LCS ranks $\phi_k(G)$. 
For example, the Chen ranks distinguish the pure braid group on 
$\ell\ge 4$ strands from the corresponding direct
product of free groups, whereas the LCS ranks don't.

As a direct application of our results, we show, in Theorem~\ref{thm:arrgts}, 
that the rational Chen Lie algebra of a complex hyperplane arrangement
is combinatorially determined.

\begin{ack}
The first author was partially supported by CERES grant 152/2003 of 
the Romanian Ministry of Education and Research.  
The second author was partially supported by NSF grants DMS-0105342 
and DMS-0311142; he is grateful to the ``Simion Stoilow" 
Institute of Mathematics of the Romanian Academy for 
hospitality during the completion of this work.
\end{ack}

\section{Malcev Lie algebras and exponential groups}
\label{sect:malcev}

Our approach to the various graded Lie algebras associated to a group $G$, 
such as $\gr(G)$ and $\gr\big(G/G^{(i)}\big)$, is based on 
the Malcev completion functor of Quillen \cite[Appendix A]{Qu}. We review
in this section Malcev  Lie algebras and their associated
exponential groups, leading to the general definition of the Malcev
completion of a group. 

\begin{definition}
\label{def:mlie}
A {\em Malcev Lie algebra} is a rational Lie algebra $L$,
together with a complete, descending $\Q$-vector space filtration,
$\{F_rL\}_{r \ge 1}$, such that:
\begin{enumerate}
\item \label{eq:ml1}
$F_1L=L$;
\item \label{eq:ml2}
$[F_rL,F_sL]\subset F_{r+s}L$, for all $r$ and $s$;
\item \label{eq:ml3}
the associated graded Lie algebra, $\gr(L)=\bigoplus _{r\ge 1}
F_rL/F_{r+1}L $,
is generated in degree~$1$.
\end{enumerate}
\end{definition}

Completeness of the filtration means that the topology on $L$ induced by
$\{F_rL\}_{r \ge 1}$ is Hausdorff, and that every Cauchy sequence
converges.  In other words, the canonical map to the inverse limit,
$\pi\colon L\to \varprojlim_r L/F_r L$, is a vector space isomorphism.  

For example, any nilpotent Lie algebra $L$, with lower central series
filtration $\{\Gamma_r L\}_{r \ge 1}$, is a Malcev Lie algebra.

\begin{definition}
\label{def:expgroup}
Let $L$ be a Malcev Lie algebra, with filtration $\{F_rL\}_{r \ge 1}$.
The {\em exponential group} associated to $L$
is the filtered group $\exp(L)$, with underlying set $L$,
group multiplication given by the Campbell-Hausdorff formula
\begin{equation}
\label{eq:chmult}
x \cdot y=x+y+\tfrac{1}{2}[x,y]+
\tfrac{1}{12}[x,[x,y]]+\tfrac{1}{12}[y,[y,x]]+\cdots ,
\quad \text{for $x,y\in L$},
\end{equation}
and filtration provided by the normal subgroups $\{\exp(F_rL)\}_{r \ge
1}$.
\end{definition}
The convergence of the series \eqref{eq:chmult} follows 
from condition \eqref{eq:ml2} above, together with the 
completeness of the filtration topology.

If $\rho \colon G  \to \exp (L)$ is a group homomorphism, 
then, as shown by  Lazard \cite{La}, 
 $\rho$ induces a graded Lie algebra map,
\begin{equation}
\label{lmal=eq}
\gr (\rho)\colon \gr (G)\otimes \Q \longrightarrow \gr (\exp(L))=\gr(L). 
\end{equation}  

\begin{definition}
\label{def:malcomp}
A {\em Malcev completion} of a group $G$ is
a group homomorphism, $\rho\colon G \to \exp (L)$,
where $L$ is a Malcev Lie algebra, such that,
for each $r\ge 1$, the factorization
$\rho_r\colon G/\Gamma_r G \to \exp (L/F_r L)$
has the property that:
given a homomorphism $f_r\colon G/\Gamma_r
G\to\exp(\widetilde{L})$,  where $\widetilde{L}$ is a nilpotent Malcev Lie
algebra, there is a unique lift
$\tilde{f}_r$  to $\exp (L/F_r L)$.
\end{definition}

The above universality property of Malcev completion is embodied
in the commuting diagram
\begin{equation}
\label{eq:univ}
\xymatrix{%
G \ar[r]^(.42){\rho}\ar@{>>}[d] &\exp(L)\ar@{>>}[d]\\
G/\Gamma_r G \ar[r]^(.42){\rho_r} \ar[dr]_(.45){f_r} &
\exp (L/F_r L) \ar@{..>}[d]^(.45){\tilde{f}_r} \\
& \exp (\widetilde{L})
}
\end{equation}
The Malcev Lie algebra $L$ is uniquely determined by the 
above universality property. The canonical Malcev completion, 
$\widehat{G}$, can be constructed as follows (see Quillen \cite{Qu68}, \cite{Qu}).

The group algebra $\Q G$ has a natural Hopf algebra structure,
with comultiplication given by $\Delta(g)=g\otimes g$, and
counit the augmentation map. Let $I$ be the augmentation ideal,
and let $\widehat{\Q G}=\varprojlim_r \Q G/I^r$ be the completion
of $\Q G$ with respect to the $I$-adic filtration. The Hopf algebra
structure extends to the completion.
Moreover, the Lie algebra of primitive elements in $\widehat{\Q G}$,
endowed with the inverse limit filtration, is a Malcev Lie algebra,
denoted by $M_G$. Set:
\begin{equation}
\label{eq:prim}
\widehat{G}=\exp (M_G)\, .
\end{equation}

For example, if $G$ is a nilpotent group, then $M_G$ is a nilpotent Lie algebra,
and so $\widehat{G}=G\otimes \Q$, the classical Malcev completion of $G$.
In general, $\widehat{G}=\varprojlim_n ( (G/\Gamma_n G) \otimes \Q )$.

For an arbitrary group $G$, Quillen also constructs a natural group 
homomorphism,
\begin{equation}
\label{qmal=eq}
\rho_G\colon G\longrightarrow \widehat{G} =\exp (M_G),
\end{equation} 
called the {\em functorial} Malcev completion homomorphism, 
which satisfies the universality property~\eqref{eq:univ}.  
Moreover, 
\begin{equation}
\label{grmal=eq}
\gr (\rho_G)\colon \gr (G)\otimes \Q \stackrel{\sim}{\longrightarrow} \gr (M_G)
\end{equation}
is an isomorphism of graded Lie algebras.  

Now, let $L$ be a Malcev Lie algebra as in Definition~\ref{def:malcomp}. 
Then, by the above discussion,  $L\cong M_G$. Hence, 
by \eqref{grmal=eq},  
\begin{equation}
\label{grmalc}
\gr (L)\cong \gr (G)\otimes \Q.
\end{equation}

\section{Malcev completion and derived series}
\label{sect:dermalcev}

We now investigate the relationship between Malcev completion
and derived series. We will need a few lemmas.

\begin{lemma}
\label{lem:m1}
Let $\rho\colon G \to \exp (L)$ be a group homomorphism.
Then $\rho\big(G^{(i)}\big)\subset \exp \big( \overline{L^{(i)}}\big)$,
for all $i\ge 1$, where $\overline{L^{(i)}}$ is the closure of $L^{(i)}$ 
with respect to the filtration topology on $L$.

\end{lemma}

\begin{proof}
Suppose $\rho(g)=x$ and $\rho(h)=y$.  By the Campbell-Hausdorff formula,
we have:
\[
\rho(ghg^{-1}h^{-1})=xyx^{-1}y^{-1}=[x,y]+
\cdots \in \exp(\overline{L'}).
\]
Thus, $\rho(G')\subset \exp (\overline{L'})$.
The general case follows by induction on $i$.
\end{proof}

In particular, the above Lemma implies that $\rho$ factors 
through a homomorphism
\[
\rho^{(i)} \colon G/G^{(i)} \longrightarrow \exp (L/\overline{L^{(i)}})\, ,
\]
for all $i\ge 1$. Here, $L/\overline{L^{(i)}}$ is endowed with the filtration 
induced from $L$; this is a Malcev filtration (in the sense of 
Definition~\ref{def:mlie}), since $\overline{L^{(i)}}$ is a closed ideal of $L$.

\begin{lemma}
\label{lem:m2}
Let $L$ be a nilpotent Malcev Lie algebra.
Then $\exp (L)^{(i)} = \exp \big(L^{(i)}\big)$, for all $i\ge 1$.
\end{lemma}

\begin{proof}
It is enough to prove the equality for $i=1$;
the general case follows by induction on $i$.
To prove the inclusion $\exp (L)' \subset \exp (L')$,
use the Campbell-Hausdorff formula
to express the group commutator $(x,y)$ as a sum
of Lie brackets, as above.
To prove the reverse inclusion, use the Zassenhaus formula
(cf.~\cite{La})
to express the Lie bracket $[x,y]$ as a product
of group commutators.
\end{proof}

We will also need the following criterion of Quillen \cite{Qu},
which characterizes the Malcev completion of a nilpotent group.
In the statement below, one has to take into account the fact that
the exponential construction described in Section~\ref{sect:malcev}
establishes a categorical equivalence between nilpotent $\Q$--Lie
algebras, and nilpotent uniquely divisible groups; see~\cite{Qu}.

\begin{prop}[Quillen]
\label{Q-crit}
Let $N$ be a nilpotent group.  A homomorphism $\varphi\colon N\to M$
defines a Malcev completion of $N$ if and only if:
\begin{enumerate}
\item $M$ is a uniquely divisible, nilpotent group;
\item $\ker \varphi$ is a torsion group;
\item for all $y\in M$, there is an integer $n\ne 0$ such that
$y^n\in \im \varphi$.
\end{enumerate}
\end{prop}

\begin{lemma}
\label{lem:m3}
Let $N$ be a nilpotent group.  Suppose $\varphi\colon N\to M$
is a Malcev completion homomorphism.  Then the restrictions to derived
series
terms, $\varphi^{(i)}\colon N^{(i)}\to M^{(i)}$, are also Malcev
completion homomorphisms.
\end{lemma}

\begin{proof}
By \cite{HMR}, the restriction of $\varphi$ to lower central series
terms,
$\Gamma_k\varphi\colon \Gamma_k N\to \Gamma_k M$, is a Malcev completion
homomorphism,
for each $k\ge 1$.  Since $\varphi^{(1)}=\Gamma_2 \varphi$, the claim is
proved for $i=1$.  The general case follows by induction on $i$.
\end{proof}

\begin{theorem}
\label{thm:Malcompl}
Let $G$ be a group, with Malcev completion $\rho \colon G \to \exp (L)$, as
in {\em Definition~\ref{def:malcomp}.}  Then, the Malcev completion of the 
derived quotient $G/G^{(i)}$ is given by
\begin{equation}
\label{eq:malcevk}
    \widehat{G/G^{(i)}} = \exp \big(L/\overline{L^{(i)}}\big), \quad 
\text{for all $i\ge 1$} .
\end{equation}
\end{theorem}

\begin{proof}
Fix $i\ge 1$. By Lemma~\ref{lem:m1}, $\rho$ factors through
a homomorphism $\rho^{(i)}\colon G/G^{(i)}\to \exp (L/\overline{L^{(i)}})$.
It remains to verify that this homomorphism satisfies the universality
property \eqref{eq:univ}.

So fix $r\ge 1$, and consider the following diagram:
\begin{equation*}
\xymatrixrowsep{10pt}
\xymatrixcolsep{-50pt}
\xymatrix@!C{
& G\ar[rrr]^{\rho} \ar[ddl]^{} \ar'[dd][ddd] &&
& \exp (L)  \ar[ddl]^{}\ar[ddd]
\\
&&&&
\\
G/G^{(i)} \ar[rrr]^{\rho^{(i)}}\ar[ddd]
&&
& \exp \left(L/\overline{L^{(i)}}\right) \ar[ddd]
\\
& G/\Gamma_r G \ar[ddl] \ar'[rr]^(.75){\rho_r}[rrr]
\ar[dr]^{f_r}&&
& \exp \left(L/F_r L\right)\ar[ddl]^{}  \ar@{-->}[dll]_(.63){\tilde{f}_r}
|!{[d];[ll]} \hole
\\
&&\exp\big(\widetilde{L}\big)&&
\\
G/\Gamma_r G \cdot G^{(i)} \ar[urr]^(.55){f} \ar[rrr]
\ar[rrr]^{\rho^{(i)}_r}&&
& \exp \left(L/(\overline{L^{(i)}}+F_r L)\right)  \ar@{..>}[ul]_(.5){\tilde{f}}
}
\end{equation*}

The maps $\rho^{(i)}$ and $\rho_r$ factor through a common quotient map,
$\rho^{(i)}_r$.  We are given a map $f$ and must lift it to $\tilde{f}$,
as indicated in the diagram.  Let $f_r$ be the composite of $f$ with
the projection map from $G/\Gamma_r G$.
The map $f_r$ has a unique lift $\tilde{f}_r$, by the
universality property of $\rho$. To construct the lift $\tilde{f}$, it is
enough to show that $(L^{(i)}+ F_r L)/F_r L$ is contained in the kernel
of $\tilde{f}_r$. (Since the topology of $L/F_r L$ is discrete,
 $(\overline{L^{(i)}}+ F_r L)/F_r L =  (L^{(i)}+ F_r L)/F_r L$.)

By Lemmas~\ref{lem:m2} and \ref{lem:m3}, we have the following:
For each $x\in (L^{(i)}+F_r L)/F_r L$, there is a
 $g\in (G^{(i)}\cdot \Gamma_r G)/\Gamma_r G$
such that $x^n=\rho_r(g)$, for some $n\ne 0$. Hence:
\[
\tilde{f}_r(x)^n=\tilde{f}_r(\rho_r(g))=f_r(g)=f(g)=1 \in
\exp(\widetilde{L})
\]
and so $\tilde{f}_r(x)=1$, since $\exp(\widetilde{L})$ is uniquely
divisible. 

To show that the lift $\tilde{f}$ is unique, let $x\in \exp(L/F_r L)$.
As before, write $x^m=\rho_r(h)$, for some $m\ne 0$ and 
$h\in G/\Gamma_r G$.
Then $\tilde{f}(x)^m=f(h)$,  and the uniqueness of $\tilde{f}$
follows from the unique divisibility  of $\exp(\widetilde{L})$.
\end{proof}

\section{Holonomy Lie algebras and formality}
\label{sect:holoformal}

In this section, we establish our main result.   
Under a formality assumption, we identify the derived quotients 
of the rational holonomy Lie algebra of a group $G$ with the associated
graded rational Lie algebra of the derived quotients of the group. 

We will only consider spaces $X$ having the homotopy  
type of a connected $CW$-complex with finite $2$-skeleton.  
Accordingly, if $G$ is the fundamental group of $X$, then  
$G$ is finitely-presented, and the classifying space $K(G, 1)$ 
can be chosen to have finite $2$-skeleton.

Let $\HH(G;\Q)$ be Chen's {\em rational holonomy Lie algebra} of
the space $K(G, 1)$,
\begin{equation}
\label{qholg=eq}
\HH(G;\Q) := \LL^*(H_1(G; \Q))/\, \text{ideal} \, \big(\im (\partial^\Q_G)\big).
\end{equation}
Here $\LL^*$ denotes the free, graded Lie algebra functor,
where the grading is given by bracket length, and 
$\partial^\Q_G\colon H_2(G;\Q)\to \LL^2(H_1(G;\Q))$ 
is the dual of the rational cup-product map, 
$\cup_G^{\Q}\colon H^1(G;\Q)\wedge H^1(G;\Q) \to H^2(G;\Q)$, with
the standard identification of $\bigwedge^2 H_1$ with $\LL^2(H_1)$.

The next definition will be crucial for our purposes. 

\begin{definition}
\label{def:1formal}
A group $G$ is called {\em $1$-formal}\/
if the Malcev Lie algebra of $G$ is isomorphic to the rational
holonomy Lie algebra of $G$, completed with respect to bracket length.
\end{definition}

In other words, $G$ is $1$-formal if 
\begin{equation}
\label{eq:1fdef}
\widehat{G}= \exp\left( \widehat{\HH(G;\Q)} \right).
\end{equation}

We are now in position to state and prove our main result (Theorem~\ref{thm1} 
from the Introduction).

\begin{theorem}
\label{thm:hatformal}
Let $G$ be a $1$-formal group. Then, for all $i\ge 1$,
\begin{equation}
\label{eq:1fiso}
\gr (G/G^{(i)})\otimes \Q =  \HH(G;\Q) / \HH(G;\Q)^{(i)}.
\end{equation}
\end{theorem}

\begin{proof}
Set $\HH = \HH(G;\Q)$ and $L=\widehat{\HH}$.  
By $1$-formality, $\widehat{G}=\exp(L)$.
By Theorem~\ref{thm:Malcompl}, we have
$\widehat{G/G^{(i)}} = \exp (L/\overline{L^{(i)}})$. 
Using the isomorphism \eqref{grmalc}, we find 
\begin{equation}
\label{iso1}
\gr (G/G^{(i)})\otimes \Q \cong \gr \big(L/\overline{L^{(i)}}\big).
\end{equation}

Now consider the canonical Lie algebra map,
$\iota \colon \HH \to L$, which sends the bracket length filtration of $\HH$
to the Malcev filtration of $L$, and embeds $\HH$ into its completion. 
For each $i\ge 1$, there is an induced map of filtered Lie algebras,
$\iota^{(i)} \colon \HH /\HH^{(i)} \longrightarrow L/\overline{L^{(i)}}$. 
It is straightforward to check that 
\begin{equation}
\label{iso2}
\gr \big(\iota^{(i)}\big) \colon \gr  \big(\HH/\HH^{(i)}\big) 
\longrightarrow \gr \big(L/\overline{L^{(i)}}\big) 
\end{equation}
is an isomorphism. Moreover, $\gr (\HH/\HH^{(i)})\cong \HH/\HH^{(i)}$, since the
filtration of $\HH/\HH^{(i)}$ is induced by the grading. Combining this isomorphism  
with those from \eqref{iso1} and \eqref{iso2} completes the proof.
\end{proof}

Here is a quick application of our result.

\begin{corollary}
\label{lhilb=cor}
Let $\LL_n$ be the free $\Z$-Lie algebra on $n$ generators (in degree $1$). 
Then $\LL_n/\LL''_n$ is a graded, free abelian group, with Hilbert series
\begin{equation}
\label{eq:chenfree}
\Hilb(\LL_n/\LL''_n,t)=nt + \sum_{k\ge 2} (k-1)\binom{n+k-2}{k} t^k .
\end{equation}
\end{corollary}

\begin{proof}
The torsion-freeness of $\LL_n/\LL''_n$ follows from the fact that
the derived subalgebra of a free $\Z$-Lie algebra is again free;
see~\cite[Exercise 11, p.77]{Bb}.

Let $\FF_n= \pi_1(\bigvee_n S^1)$ be the free group on $n$ generators.
Since wedges of circles are formal spaces, Theorem~\ref{thm:hatformal}
applies to $\FF_n$. Note that $\HH (\FF_n ;\Q)= \LL_n \otimes \Q$, 
as follows directly from definition~\eqref{qholg=eq}. 
Therefore, the coefficients of the Hilbert series of $\LL_n/\LL''_n$ are
equal to the Chen ranks, $\theta_k (\FF_n)$.  Those ranks were computed  
by Chen \cite{Ch51}  (see also \cite{Mu, MT, CS2}).  Putting things together 
yields \eqref{eq:chenfree}.
\end{proof}

Theorem~\ref{thm:hatformal} can be used as 
a $1$-formality test.  We illustrate with a simple example. 
Another example (for link groups) will be given in \S\ref{subsect:chenmura}. 

\begin{example}
\label{notf=ex}
Let $G =\langle x_1, x_2 \mid ((x_1,x_2),x_2) \rangle$. 
Plainly, $\theta_3 =\phi_3 \le 1$. On the other hand, 
$\HH(G;\Q) =\LL_2\otimes \Q$ (this follows from Proposition~\ref{mod3=prop}).  
By Corollary~\ref{lhilb=cor}, the rank of the degree $3$ piece of 
$\HH(G;\Q)/\HH(G;\Q)''$ is $2$. We infer that $G$ cannot be $1$-formal, 
since this would contradict Theorem~\ref{thm:hatformal}.
\end{example}

\section{Integral holonomy Lie algebras}
\label{zhol=sec}

Our next objective is to find an integral version of Theorem~\ref{thm:hatformal}. 
We start by recalling the definition and basic properties of 
$\Z$-holonomy Lie algebras.  As before, let $X$ be a space having the 
homotopy type of a connected CW-complex with finite $2$-skeleton. 

\begin{definition}[\cite{MP}]
\label{holz=def}
Assume $H_1 X=H_1(X; \Z)$ is torsion-free. 
The {\em integral holonomy Lie algebra} of $X$ is
\[
\HH(X) := \LL^*(H_1 X)/\, \text{ideal}\, (\im \left(\overline{\partial}_X)\right).
\]
Here $\overline{\partial}_X := \partial_X \circ \kappa$, where
$\partial_X \colon \Hom (H^2 X, \Z)\to \bigwedge^2 H_1 X = \LL^2(H_1 X)$
is dual to the cup-product map $\cup_X \colon \bigwedge^2 H^1 X \to H^2 X$, and
$\kappa \colon H_2 X \surj \Hom (H^2 X, \Z)$ is the canonical surjection induced by
the Kronecker pairing. If $G$ is a finitely-presented group with torsion-free 
abelianization, then the $\Z$-holonomy Lie algebra of $G$ is
\[
\HH(G):= \HH(K(G, 1)).
\] 
\end{definition}

From the above definition, we see that $\HH(X)$ is a $\Z$-form of
$\HH (X; \Q)$, that is,
\begin{equation}
\label{zform=eq}
\HH(X; \Q) = \HH(X) \otimes \Q,
\end{equation}
as graded Lie algebras. Obviously, both $\overline{\partial}_X$ and
$\partial^{\Q}_X$ are natural with respect to continuous maps.
Consider now a classifying map, $f\colon X \to K(G, 1)$, where
$G = \pi_1 (X)$. Since $H_1 f$ is an isomorphism and $H_2 f$ is an
epimorphism, we infer that $f$ induces an isomorphism between 
$\Z$-holonomy Lie algebras:
\begin{equation}
\label{id=eq}
\HH(X)\cong \HH(\pi_1 (X)). 
\end{equation}
Similarly, $\HH(X;\Q)\cong \HH(\pi_1(X);\Q)$.  

By the universality property of free Lie algebras, the standard identification 
$H_1 G = \gr_1 (G)$ extends to a surjective map of graded $\Z$-Lie algebras,
\begin{equation*}
\Psi_G \colon \LL^*(H_1 G) \surj \gr (G),
\end{equation*}
natural with respect to group homomorphisms. 

\begin{prop}[\cite{MP}]
\label{jm=prop}
The above map $\Psi_G$ sends $\im (\overline{\partial}_G) = \im (\partial_G)$
(see {\em Definition~\ref{holz=def}}) to zero, thus inducing a natural epimorphism of 
graded $\Z$-Lie algebras,
\begin{equation}
\label{natgr=eq}
\Psi_G \colon \HH(G) \surj \gr (G)\, ,
\end{equation}
for any finitely-presented group $G$ with torsion-free abelianization.
\end{prop}

We now show that $\Psi_G$ further descends to the derived quotients.

\begin{prop}
\label{prop:nat}
For each $i\ge 1$, there is a natural epimorphism of graded $\Z$-Lie algebras,
$\Psi_G^{(i)}\colon \HH(G)/\HH(G)^{(i)} \surj \gr(G/G^{(i)})$,
which fits into the commuting diagram
\begin{equation*}
\xymatrix{%
\HH(G) \ar@{>>}[rr]^(.5){\Psi_G} \ar@{>>}[d] & & \gr(G)
\ar@{>>}[d]^{}  \\
\HH(G)/\HH(G)^{(i)} \ar@{>>}[rr]^(.5){\Psi_G^{(i)} }& &
\gr\!\big(G/G^{(i)}\big)
}
\end{equation*}
\end{prop}

\begin{proof}
Clearly, $\Psi_G\big( \HH(G)^{(i)} \big) \subset \gr(G)^{(i)}$.
Furthermore, the canonical projection $\gr(G)\surj \gr(G/G^{(i)})$
takes $\gr(G)^{(i)}$ to $0$, since the Lie bracket on
$\gr (G/G^{(i)})$  is induced by the group commutator.
Thus, $\Psi_G$ passes to the quotients. The induced map,
$\Psi_G^{(i)}$, is automatically surjective.
\end{proof}

\begin{corollary}
\label{thm:1formal}
Let $G$ be a $1$-formal group, with torsion-free abelianization.
Then the natural Lie algebra maps
\[
\Psi_G^{(i)}\otimes \Q \colon
\left(\HH(G)/\HH(G)^{(i)}\right)\otimes \Q \longrightarrow
\gr\!\big(G/G^{(i)}\big)\otimes \Q
\]
are isomorphisms, for all $i\ge 1$.
\end{corollary}

\begin{proof}
By Theorem~\ref{thm:hatformal}, the graded Lie algebras
$\gr (G/G^{(i)})\otimes \Q$ and
$\big(\HH(G)/\HH(G)^{(i)}\big)\otimes \Q$ are isomorphic.  Therefore, the
linear map $\Psi_G^{(i)}\otimes \Q $ induces in each degree 
a surjection between
$\Q$-vector spaces of the same (finite) dimension,
hence an isomorphism.
\end{proof}

The following Proposition, together with Corollary~\ref{thm:1formal}, 
implies Theorem~\ref{thm2} from the Introduction.

\begin{prop}
\label{tors}
Suppose  $\Psi_G^{(i)}\otimes \Q$ is an isomorphism.
Then $\Psi_G^{(i)}$ induces a surjection on $p$-torsion, for all primes
$p$.  If, moreover, the Hilbert series of the graded vector space
$\big( \HH(G)/\HH(G)^{(i)}\big) \otimes \F_p$ is independent of $p$,
then $\Psi_G^{(i)}$ is an isomorphism of graded $\Z$-Lie algebras,
and the graded abelian group $\gr\big(G/G^{(i)}\big)$ is torsion-free.
\end{prop}

Of course, the same statement also holds for
$\Psi_G \colon \HH(G) \surj \gr (G)$; applications may be found 
in~\cite{BP1, MP}.  Here is our first application for $\Psi_G^{(i)}$.

\begin{corollary}
\label{hichen=cor}
For each $i\ge 1$, there is an isomorphism 
$\gr \big(\FF_n/\FF_n^{(i)}\big)\cong \LL_n/\LL_n^{(i)}$ 
of graded $\Z$-Lie algebras.  Moreover, $\gr \big(\FF_n/\FF_n^{(i)}\big)$
is torsion-free, as a graded abelian group.	
\end{corollary}

\begin{proof}
As in the proof of Corollary~\ref{lhilb=cor}, start by noting that
$\FF_n$ is $1$-formal, $\HH(\FF_n)= \LL_n$, and $\LL_n/\LL_n^{(i)}$ is
torsion-free. Then apply Proposition~\ref{tors}.
\end{proof}  

\section{The infinitesimal Alexander invariant of a quadratic Lie
algebra}
\label{sect:holoalex}

Let $\LL_Y$ be the the free Lie algebra (over $\Z$) 
on a set $Y$.  Denote by $S=\Z[Y]$ the symmetric algebra on $Y$, 
with variables in degree $1$. By the Poincar\'{e}-Birkhoff-Witt Theorem, 
$S$ is isomorphic to the universal enveloping algebra $U(\LL_Y/\LL'_Y)$. 

\begin{definition}
\label{def:alex}
Let $E$ be a graded Lie algebra.  Suppose $E=\LL_Y/I$,
where $I$ is a homogeneous ideal, generated in
degree $\ge 2$. The {\em infinitesimal Alexander invariant} of
$E$ is the graded $S$-module $B(E) :=E'/E''$, with module structure given
by
the adjoint representation, via the exact sequence of Lie algebras
\begin{equation}
\label{eq:liealex}
0 \to E'/E'' \to E/E'' \to E/E' \to 0.
\end{equation}
\end{definition}

In the case when $Y$ is finite and $I$ is generated in degree $2$,
we can write down explicitly a finite presentation for $B(E)$.
We shall use the Koszul resolution
$(S\otimes \bigwedge^k_Y, d_k)$ of the trivial $S$-module $\Z$, where
$\bigwedge^*_Y$ denotes the exterior algebra on $Y$.
We will also identify $\LL_Y^1=\bigwedge_Y^1$ and
$\LL_Y^2=\bigwedge_Y^2$,  with $[x,y]$ corresponding to $x\wedge y$.

\begin{theorem}
\label{thm:ppres}
Let $E$ be a graded Lie algebra.  Suppose $E$ admits a quadratic,
finite presentation of the form $E=\LL/I$,
where $\LL$ is the free Lie algebra on a finite set $Y$
and $I$ is a graded ideal generated by  $J := I \cap \LL^2$.
Then
\begin{equation}
\label{eq:delta}
S\otimes \left(\bigwedge\nolimits^3_Y \oplus J\right) \xrightarrow{
  d_3 + (\id\otimes \partial) } S\otimes \bigwedge\nolimits^2_Y
\longrightarrow B(E) \to 0
\end{equation}
is a homogeneous, finite $S$-presentation of $B(E)$, where:
\begin{itemize}
\item $\deg (\bigwedge^k_Y)=k$, $\deg(J)=2$,
\item $\partial$ is the inclusion of $J$ into $\bigwedge_Y^2$,
\item $d_3 (x\wedge y\wedge z)= x\otimes y\wedge z -
y\otimes x\wedge z + z\otimes x\wedge y$, for $x,y,z \in Y$.
\end{itemize}
\end{theorem}

\begin{proof}
We start by finding a presentation for the infinitesimal
Alexander invariant of the free Lie algebra $\LL=\LL_Y$
on the set $Y=\{y_1,\dots,y_n\}$.  Define an $S$-linear map
\begin{equation}
\label{eq:eta}
\eta \colon S\otimes \bigwedge\nolimits^2_Y \longrightarrow \LL'/\LL'' ,
\qquad  x\wedge y \mapsto [x,y]\bmod \LL''.
\end{equation}

We claim that $\eta$ is surjective.  Indeed,  let $[u,v]\in\LL'$.
If both $u$ and $v$ have bracket length greater than $1$,
then $[u,v]\in \LL''$. So we may assume that $u$
has length $1$. Now, if $v$ also has length $1$,        
then $[u,v]=\eta (u\wedge v)$. The claim follows by induction
on the bracket length of $v$, using commutator calculus.

Now let $d_3\colon S\otimes \bigwedge\nolimits^3_Y  \to
S\otimes \bigwedge\nolimits^2_Y$ be the Koszul differential.
By the Jacobi identity, we have $\eta\circ d_3=0$.
 From Corollary~\ref{lhilb=cor}, we know that the degree $k$ component
of $\LL'/\LL''$ is free abelian, of rank equal to
$\theta_k(\FF_n)=(k-1) \binom{n+k-2}{k}$,
for each $k\ge 2$.   A standard dimension-counting argument
shows that this rank is the same as the rank of the degree $k$ component
of $\coker d_3$, which is also free abelian. 
Thus, $\ker \eta=\im d_3$ and the following sequence is exact:
\begin{equation}
\label{eq:blp}
S\otimes \bigwedge\nolimits^3_Y  \xrightarrow{\ d_3\ }
S\otimes \bigwedge\nolimits^2_Y \xrightarrow{\ \eta\ } B(\LL_Y) \to 0.
\end{equation}

We return now to the proof in the general case: $E=\LL/I$,
where $I$ is generated in degree two. Let $\pi\colon \LL\to E$ 
be the canonical projection.  Note that $E'/E''=\LL'/(\LL''+I)$, and so
we have an exact sequence
\begin{equation}
\label{eq:exact}
0\to (I+\LL'')/\LL'' \to \LL'/\LL'' \xrightarrow{\ \pi \ }
E'/E''\to 0.
\end{equation}
It is readily seen that the $S$-linear surjection $\pi\circ \eta\colon
S\otimes \bigwedge\nolimits^2_Y \to E'/E''$ has kernel
equal to the image of $d_3 + (\id\otimes \partial)\colon
S\otimes \big(\bigwedge\nolimits^3_Y \oplus J\big) \to
S\otimes\bigwedge\nolimits^2_Y$. This finishes the proof.
\end{proof}

From the Hilbert-Serre Theorem, we also obtain the following corollary.

\begin{corollary}
\label{cor:ratl}
If $E$ is a quadratic, finitely presented, graded Lie algebra over
a field $\K$, then
\begin{equation}
\label{eq:addhilb}
\Hilb_\K (E/E'',t)=\Hilb_\K (E'/E'',t)+\Hilb_\K (E/E',t)
\end{equation}
is a {\em rational function} in $t$.
\end{corollary}

The main example to keep in mind is that of the holonomy Lie algebra,
$E=\HH(G)$, of a finitely presented group $G$
with torsion-free abelianization. 
We simply call $B(\HH(G))=\HH(G)'/\HH(G)''$ the
{\em infinitesimal Alexander invariant} of $G$.

\section{Chen Lie algebras of one-relator groups}
\label{sec:onerel}

In this section, we isolate a class of groups $G$ for which the holonomy Lie 
algebra $\HH(G)$ may be described directly in terms of the defining 
relations. If $G$ has a single relator (of a certain type), we determine 
the Chen Lie algebra $\gr(G/G'')$. 

\begin{definition}
\label{def:commrel}
A group $G$ is called a {\em commutator-relators} group if
it admits a finite presentation of the form
$G=\langle x_1,\dots ,x_n\mid r_1,\dots ,r_m\rangle$,
where each relator $r_j$ belongs to the commutator subgroup
of the free group $\FF_n=\langle x_1,\dots ,x_n\rangle$.
\end{definition}

Note that $G/G'=\Z^n$, with canonical $\Z$-basis, $Y=\{ y_1, \ldots, y_n \}$,
given by the classes of the group generators. Recall that
$\gr (\FF_n)= \LL_Y$. Denote by $\{ \r_j \in \LL_Y^2 \}_{1\le j \le m}$ 
the images of the relators, modulo $\Gamma_3 \FF_n$. 

\begin{prop} 
\label{mod3=prop}
If $G=\langle x_1,\dots ,x_n\mid r_1,\dots ,r_m\rangle$ is a
commutator-relators group, then
\begin{equation}
\label{eq:holo-one}
\HH(G) = \LL (y_1,\ldots, y_n)/\, {\rm ideal} \left(\r_1,\ldots, \r_m\right).
\end{equation}
\end{prop}

\begin{proof}
Let $\Delta_G \colon H_2 G \to H_2 (G\times G)$ be the homomorphism 
induced by the diagonal map $G\to G\times G$, and let  
$\pr \colon H_2 (G\times G) \to H_1 G \otimes H_1 G$ be the 
projection given by the K\"{u}nneth formula.  Note that 
$\Delta_G$ commutes with the automorphism of $H_2(G\times G)$ 
induced by the flip map  $G\times G \to G\times G$. Thus, the 
composite $\pr \circ \Delta_G$ takes values in $\bigwedge^2 H_1 G$, 
viewed as a subgroup of $\bigotimes^2 H_1 G$ via the embedding 
$x\wedge y\mapsto x\otimes y - y\otimes x$.  So we may define 
the reduced diagonal, 
$\overline{\Delta}_G \colon H_2 G \to H_1 G \wedge H_1 G$, 
as the co-restriction of  $\pr \circ \Delta_G$. 
It is now straightforward to check that $\overline{\Delta}_G$ 
coincides with the map $\overline{\partial}_G$ from 
Definition~\ref{holz=def}. Equality \eqref{eq:holo-one} follows 
at once from \cite[Lemma 2.8]{P}.
\end{proof}

Let $X_g$ be a closed, orientable surface of genus $g\ge 1$,
with fundamental group
$\pi_1(X_g)=\langle x_1,x'_1,\dots ,x_g,x'_g \mid
\prod_{i=1}^g (x_i,x'_i)\rangle$. In \cite{L70},  Labute 
showed that the associated graded Lie algebra of
$\pi_1(X_g)$ has presentation
\begin{equation}
\label{eq:hliesurf}
\gr(\pi_1(X_g))=\LL(x_1,x'_1,\dots ,x_g,x'_g )/([x_1,x'_1]+\cdots
+[x_{g},x'_{g}]), 
\end{equation}
with graded ranks $\phi_k=\phi_k(\pi_1(X_g))$ given by
\begin{equation}
\label{eq:lcssurf}
\prod_{k=1}^{\infty}(1-t^k)^{\phi_k}=1-2gt+t^2.
\end{equation}
(This formula can also be derived from \cite{PY}, using the fact  
that $X_g$ is  formal, and $H^*(X_g; \Q)$ is a Koszul algebra.) 

In our next result, we provide analogs of \eqref{eq:hliesurf}
and~\eqref{eq:lcssurf} above, for the Chen Lie algebras of 
``surface-like" one-relator groups.

\begin{theorem}
\label{1rel=thm}
Let $G=\langle x_1, \dots, x_n \mid r\rangle$ be a commutator-relators
group defined by a single relation. Denote by
$\r\in \LL^2 (y_1, \dots, y_n)$ the class of $r$ modulo $\Gamma_3$.
Assume that $\r \not\equiv 0$ (modulo $p$), for every prime $p$.
If $G$ is $1$-formal, then
\begin{equation}
\label{chpres=eq}
\gr (G/G'') = \LL (y_1, \dots,y_n)/ 
\big(  {\rm ideal} (\r) +\LL'' (y_1,\dots,y_n)\big),
\end{equation}
as graded $\Z$-Lie algebras.  Moreover, $\gr (G/G'')$ is torsion-free as 
a graded abelian group, and its Hilbert series, 
$\Hilb(\gr (G/G''),t)=\sum_{k\ge 1} \theta_k t^k$, is given by
\begin{equation}
\label{hsurf=eq}
\Hilb(\gr (G/G''),t)=1+ nt -\frac{1-nt+t^2}{(1-t)^n}.
\end{equation}
\end{theorem}

\begin{proof}
By Proposition~\ref{mod3=prop}, the graded $\Z$-Lie algebra on the
right-hand side of equation~\eqref{chpres=eq} equals $\HH(G)/\HH(G)''$.
By Corollary~\ref{thm:1formal}, Proposition~\ref{tors} and 
Theorem~\ref{thm:ppres}, all we have to do is to check that
\begin{equation}
\label{indp=eq}
\Hilb_{\F_p}(B(\HH(G))\otimes \F_p ,t)= 1- \frac{1-nt+t^2}{(1-t)^n}\, ,
\end{equation}
for all primes $p$. Fix then $\F_p$-coefficients (omitting them from notation).

From presentation~\eqref{eq:delta} and the exactness of Koszul resolutions,
we infer that
\begin{equation*} 
\label{mupres=eq}
B(\HH(G)) =\coker \, \{ S\xrightarrow{\mu} \ker \, (d_1)\} ,
\end{equation*}
where $\mu$ denotes multiplication by $d_2 (1\otimes \r)\in \ker (d_1)$.
Write $1\otimes \r =\sum_{1\le i<j\le n} c_{ij}
\otimes y_i \wedge y_j$, with $c_{ij}\in \F_p$. Then
\[
d_2 (1\otimes \r) = \sum\limits_{i<j} c_{ij}
(y_i\otimes y_j -y_j\otimes y_i) ,
\]
and this is non-zero, since $\r \neq 0$ by assumption. Hence $\mu$ 
is injective, and formula \eqref{indp=eq} readily follows.
\end{proof}

For many one-relator groups $G$, the above theorem 
determines the structure of the Chen Lie algebra $\gr(G/G'')$, 
and computes the Chen ranks $\theta_k(G)$.  
We illustrate with several examples. 

\begin{example}
\label{curve=ex}
Let $G$ be the fundamental group of an irreducible curve in the 
complex projective plane.   
Denote by $\Sigma \subset X$ the set of singular points of $X$. 
For each $x\in \Sigma$, let $m_x$ be the number of local branches 
passing through $x$.
Set $m:=\sum_{x\in \Sigma} (m_x -1)$, and $n:= 2g+m$, where 
$g$ is the genus of $X$. Then, as is well-known, $X$ has the 
homotopy type of the formal space $X_g \vee (\bigvee_m S^1)$. 
If $g=0$, then $G=\FF_m$, and Corollaries~\ref{lhilb=cor} and \ref{hichen=cor} 
apply.  If $g\ge 1$, then $G=\pi_1(X_g)*\FF_m$ has a presentation with $n$
generators and one commutator relator, and Theorem~\ref{1rel=thm} applies. 
\end{example}

\begin{example}
\label{psurf=ex}
Let $G=\langle x_1, \dots, x_n \mid r\rangle$ be a one-commutator-relator group. 
Assume that $\r\in \bigwedge^2 (y_1,\dots, y_n)$ is non-degenerate, and 
$\r \not\equiv 0$ ($\bmod$ $p$), for every prime $p$.   
Denote by $X$ the finite $2$-complex 
associated to the given presentation of $G$. By non-degeneracy, 
$n=2g$ and $H^*(X; \Q)=H^*(X_g; \Q)$, as graded algebras. By the
rigidity results from~\cite[\S2.5]{BP2}, the group $G=\pi_1(X)$ is $1$-formal.  
Hence, Theorem~\ref{1rel=thm} applies.
\end{example}

\section{Alexander invariant}
\label{sect:alexinv}

We now turn to the relationship between the holonomy Lie algebra 
and the Alexander invariant. 
We start by reviewing some basic material on the Alexander module
and the Alexander invariant of a group, as well as their associated
graded modules, and their connection to the Chen groups.

\subsection{Alexander modules}
\label{subsec:ab}
Let $G$ be a finitely-presented group.  Let $\Z G$ be the group-ring,
$\epsilon\colon \Z G\to \Z$ the augmentation map, given by
$\epsilon (\sum n_g g)=\sum n_g$, and $I_G:=\ker \epsilon$ the
augmentation ideal. Finally,
let $G/G'$ be the abelianization of $G$, and $\alpha\colon G\surj G/G'$
the canonical projection.

Associated to $G$ there are two important modules over $\Z(G/G')$:
\begin{enumerate}
\item The {\em Alexander module}: $A_G=\Z(G/G') \otimes_{\Z G} I_G$,
the module induced from $I_G$ by the extension of $\alpha$ to
group-rings.

\item The {\em Alexander invariant}: $B_G=G'/G''$, with $G/G'$ acting
on the cosets of $G''$ via conjugation: $gG'\cdot hG''=ghg^{-1}G''$,
for $g\in G$, $h\in G'$.
\end{enumerate}
These two $\Z(G/G')$-modules fit into the Crowell exact sequence,
\begin{equation}
\label{eq:crowell}
0\to B_G\to A_G\to I_{G/G'}\to 0.
\end{equation}

If $X$ is a connected CW-complex with finite $2$-skeleton 
and $G=\pi_{1}(X,*)$, and if
$p\colon \overline{X}\to X$ is the maximal abelian cover, then
the homology exact sequence of the pair
$(\overline{X}, p^{-1}(*))$ splits off the sequence \eqref{eq:crowell},
with $B_G=H_{1}(\overline{X})$ and $A_G=H_{1}(\overline{X}, p^{-1}(*))$.

\subsection{Associated graded modules}
\label{subsec:grB}
Set $I:=I_{G/G'}$.
The module $B_G$ comes endowed with the $I$-adic filtration,
$\{I^k B_G\}_{k\ge 0}$.  Let $\gr( B_G)=\bigoplus _{k\ge 0} I^{k}B_G/
I^{k+1}B_G$
be the associated graded module over the ring
$\gr ( \Z (G/G'))=\bigoplus _{k\ge 0} I^{k}/ I^{k+1}$.
Then, as shown by W.S.~Massey \cite{Ms},
\begin{equation}
\label{eq:massey}
\gr_k (G/G'')= \gr_{k-2} (B_G)
\end{equation}
for all $k\ge 2$, where recall that the associated graded on the left
side is
taken with respect to the lower central series filtration. In particular,
\begin{equation}
\label{eq:hilbmassey}
\sum_{k\ge 0} \theta_{k+2}(G)\cdot t^k = \Hilb(\gr (B_{G})\otimes \Q,t).
\end{equation}

\begin{example}
\label{ex:braid}
Let $G=P_n$ be Artin's pure braid group on $n$ strings. As shown by
Kohno \cite{Ko85}, the graded ranks $\phi_k=\phi_k(P_n)$ are given by:
\begin{equation}
\label{eq:lcspure}
\prod_{k=1}^{\infty}(1-t^k)^{\phi_k}=\prod_{j=1}^{n-1} (1-jt)  .
\end{equation}

The Hilbert series of the Alexander invariant
of $P_n$ was computed in \cite{CS1}:
\begin{equation}
\label{eq:alexpure}
\Hilb(B_{P_n},t)=
\frac{\binom{n+1}4}{(1-t)^2} - \binom{n}{4}.
\end{equation}
It follows that $\theta_k(P_n)=(k-1)\binom{n+1}{4}$, for $k\ge 3$.
\end{example}

\subsection{Presentations for Alexander modules}
\label{subsec:presalex}
Suppose $G/G'$ is torsion-free, of rank $n=b_1(G)$, and
fix a basis $\{t_1,\dots ,t_n\}$.  This identifies the
group ring $\Z(G/G')$ with the ring of Laurent polynomials
$\Lambda=\Z[t_1^{\pm 1},\dots, t_n^{\pm 1}]$.
Under this identification, the augmentation ideal $I$
corresponds to the  ideal $\I=(t_1-1,\dots, t_n-1)$.

Let $G=\langle x_1,\dots, x_n \mid r_1, \dots, r_m\rangle$ be
a commutator-relators group.
Note that $G/G'=\Z^n$. Pick as basis elements $t_i=\alpha \varphi(x_i)$,
where  $\varphi\colon\FF_n\surj G$ is the canonical projection.
The Alexander module of $G$ admits a finite presentation
\begin{equation}
\label{eq:apres}
\Lambda^m \xrightarrow{D_G=\big( \alpha\varphi \partial_i(r_j) \big)}
\Lambda^n \to A_G \to 0,
\end{equation}
where $\partial_i= \frac{\partial }{\partial x_i}\colon \FF_n\to \Z\FF_n$
are the Fox free derivatives.

Let $(\Lambda^{\binom{n}{k}},\delta_k)$ be the standard Koszul
resolution of $\Z$ over $\Lambda$.  
By the fundamental formula of Fox calculus,
$\delta_1\circ D_G=0$; see~\cite{F}. 
A diagram chase (see \cite{Ms}) gives a finite
presentation for the Alexander invariant:
\begin{equation}
\label{eq:bpres}
\Lambda^{\binom{n}{3}+m}\xrightarrow{\Delta_G=\delta_3+\vartheta_G}
\Lambda^{\binom{n}{2}} \to B_G \to 0,
\end{equation}
where $\vartheta_G$ is a map satisfying $\delta_2\circ \vartheta_G=D_G$.
(Such a map always exists, by basic homological algebra.)

\section{Linearized Alexander invariant}
\label{sect:linalex}

In this section, we show that the linearized
Alexander invariant of a commutator-relators group
coincides with the infinitesimal Alexander invariant of the
holonomy Lie algebra of the group.

\subsection{Magnus embedding}
\label{subs:magnus}

The ring $\Lambda=\Z[t_1^{\pm 1}, \dots, t_n^{\pm 1}]$
can be viewed as a subring of the formal power series ring
$P=\Z[[s_1,\dots, s_n]]$ via the Magnus embedding
$\mu\colon\Lambda \hookrightarrow P$, defined by $\mu(t_i)=1+s_i$.  
Note that $\mu$ sends the ideal $\I=(t_1-1,\dots ,t_n-1)$
to the  ideal $\m=(s_1,\dots,s_n)$.
Passing to associated graded rings (with respect to the filtrations
by powers of $\I$ and $\m$), the homomorphism $\gr(\mu)$ 
identifies $\gr(\Lambda)$ with the polynomial ring
$S=\gr(P)=\Z[s_1,\dots , s_n]$.

As a result, we may view the associated graded Alexander invariant,
$\gr (B_G)$, as an $S$-module.  To compute the Hilbert series of this 
module, one needs to know (at the very least) a finite presentation.  
But even if an explicit presentation for $B_G$ is known, finding
a presentation for $\gr( B_G)$ involves an arduous Gr\"obner basis
computation, see \cite{CS1}, \cite{CS2}. 
We turn instead to a more manageable approximation for the 
Alexander invariant of a commutator-relators group $G$.

\subsection{Linearized Alexander module}
\label{subs:linalexmod}

For each $q\ge 0$, let $\mu^{(q)}\colon \Lambda \to P/\m^{q+1}$
be the $q$-th truncation of $\mu$.
Since all the relators of $G$ are commutators, the entries of
the Alexander matrix $D_G$ from \eqref{eq:apres} are in the ideal $\I$,
and so $\mu^{(0)}(D_G)$ is the zero matrix.
Hence, all the entries of the linearized  Alexander matrix,
$D^{(1)}_G :=\mu^{(1)}(D_G)$, belong to $\m/\m^2=\gr_1(P)$,
and thus can be viewed as linear forms in the variables of $S$.

\begin{definition}
\label{def:linalexmod}
The {\em linearized Alexander module} of a commutator-relators group $G$
is the $S$-module $\Amod_G=\coker D^{(1)}_G$.
\end{definition}

The linearized Alexander matrix depends only on the relators of $G$,
modulo length~$3$ commutators.  Indeed, let $r\in \Gamma_2 \FF_n$ 
be a commutator-relator, and set 
\begin{equation*}
\label{r3}
r^{(3)}:=\prod_{i<j}(x_i, x_j)^{\epsilon_{ij}(r)},
\end{equation*}
where $\epsilon_{ij}=\epsilon \partial_i \partial_j$.  Then 
\begin{equation}
\label{rbar=eq}
r \equiv r^{(3)}, \quad \text{modulo} \quad \Gamma_3 \FF_n.
\end{equation}
 A Fox calculus computation shows that 
$\mu^{(1)}\alpha\varphi\partial_k(r)=
\mu^{(1)}\alpha\varphi\partial_k(r^{(3)})$, see~\cite{F, MS}.  
Consequently, the entries of the $m\times n$ matrix $D_G^{(1)}$ 
are given by:
\begin{equation}
\label{eq:linalex}
\left(D_G^{(1)}\right)_{k,j}=\sum_{i=1}^{n}\epsilon_{ij}(r_k) s_i .
\end{equation}

\subsection{Linearized Alexander invariant}
\label{subs:linalexinv}
As above, let  $G$ be an $n$-generator, commutator-relators group.
Identify the symmetric algebra on $H_1(G)=\Z^n$ with the polynomial
ring $S=\Z[s_1,\dots,s_n]$, and let 
$\big(S\otimes \bigwedge^k H_1(G),d_k\big)$ be the  Koszul complex 
of $S$.  The map $d_1\colon S^n \to S$ is simply the matrix of variables; 
its image is the  ideal $\m=(s_1,\dots ,s_n)$.  Notice that 
$d_1\circ D^{(1)}_G=0$.  Hence, the map $d_1$ factors through the 
cokernel of $D^{(1)}_G$, giving an epimorphism $\Amod_G\surj \m$.

\begin{definition}
\label{def:linalexinv}
The {\em linearized Alexander invariant} of a commutator-relators group
$G$ is the graded $S$-module $\B_G=\ker\big(\Amod_G\surj \m\big)$.
\end{definition}

Recall we also constructed (from the $\Z$-holonomy Lie algebra 
of $G$) a graded module over the polynomial ring $\Z[y_1,\dots, y_n]$: 
the infinitesimal Alexander invariant, $B(\HH(G))=\HH(G)'/\HH(G)''$. 
Renaming variables, $B(\HH(G))$ becomes a module over $S=\Z[s_1,\dots,s_n]$. 

\begin{prop}
\label{prop:linalex}
Let $G$ be a commutator-relators group. Then $\B_G \cong B(\HH(G))$, 
as graded $S$-modules. 
\end{prop}

\begin{proof}
By Theorem~\ref{thm:ppres}, Proposition~\ref{mod3=prop}, and exactness
of the Koszul resolution,
\begin{equation}
\label{ali1=eq}
B(\HH(G)) =\coker \left\{ d_2 \colon S\otimes J \longrightarrow \ker (d_1) \right\},
\end{equation}
where $J\subset \LL^2_Y \equiv \bigwedge^2_Y$ is generated by
$\{ \r_k \}_{1\le k \le m}$.
By definition, 
\begin{equation} 
\label{ali2=eq}
\B_G =\coker  \big\{ D_G^{(1)} \colon S^m \longrightarrow \ker (d_1) \big\} .
\end{equation}
By \eqref{eq:linalex}, \eqref{ali1=eq}, and~\eqref{ali2=eq}, 
we are left with checking that
\begin{equation}
\label{ali3=eq}
d_2 (\r) =\sum_{1\le i,j \le n}\epsilon_{ij}(r) s_i y_j, 
\end{equation}
for any $r\in \Gamma_2 \FF_n $. 
We know from \eqref{rbar=eq} that $\r=\overline{r^{(3)}}$.  Hence:
\begin{equation}
\label{ali4=eq}
\r=\sum_{1\le i<j\le n}\epsilon_{ij}(r) y_i\wedge y_j .
\end{equation}
Therefore, by the very definition of Koszul differentials, 
\begin{equation}
\label{ali5=eq}
d_2(\r)=\sum_{i<j}\epsilon_{ij}(r) (s_i y_j -s_j y_i) .
\end{equation}
Symmetry properties from the free differential calculus~\cite{F} 
show that \eqref{ali3=eq} is equivalent to 
\eqref{ali5=eq}.
\end{proof}

\begin{corollary}
\label{cor:lininfalex}
Let $G$ be a commutator-relators group,
with Alexander invariant $B_G$, linearized Alexander invariant
$\B_G$  and  infinitesimal Alexander invariant
$B(\HH(G))$. Then
\begin{equation}
\label{94=eq}
\Hilb(\B_G\otimes \Q,t)=
\Hilb(B(\HH(G))\otimes \Q,t) .
\end{equation}
If moreover $G$ is $1$-formal, then
\begin{equation}
\label{eq:hilbbg}
t^2 \Hilb(\gr( B_G)\otimes \Q,t)=
\sum_{k\ge 2} \theta_k(G) t^k=
\Hilb(B(\HH(G))\otimes \Q,t)  .
\end{equation}
\end{corollary}

\begin{proof}
The first equality from~\eqref{eq:hilbbg} follows
from Massey's isomorphism \eqref{eq:massey},
for an arbitrary finitely presented group $G$.
If $G$ is $1$-formal, the second equality is provided by
Theorem~\ref{thm:hatformal} (with $i=2$), via~\eqref{eq:liealex}.
\end{proof}

The next example shows that the Hilbert series of the linearized 
Alexander invariant does depend in general on the characteristic 
of the field over which it is computed. 

\begin{example}
\label{ex:tors}
Let $G$ be the group with generators $x_1, x_2$ and a single 
relator $r=(x_1,x_2^p)$, where $p$ is a prime. 
Since $\r=p [y_1,y_2]$ is non-degenerate, $G$ is $1$-formal, 
see Example~\ref{psurf=ex}. 
By Proposition~\ref{mod3=prop}, we have 
$\HH(G)=\LL(y_1,y_2)/\text{ideal} (p[y_1,y_2])$.  
By Theorem~\ref{thm:ppres}, the module  $\B_G=B(\HH(G))$ is 
isomorphic to $\Z[y_1,y_2]/(p)$, with degree shifted by $2$.  Hence:
\[
\Hilb(\B_G \otimes \K,t)=
\begin{cases}
0, & \text{if\: $\ch(\K)\ne p$},\\
\big(\frac{t}{1-t}\big)^2, & \text{if\: $\ch(\K)= p$}.
\end{cases}
\]
\end{example}

\section{Links in $S^3$ and the Murasugi conjecture}
\label{sect:lks}

\subsection{The holonomy Lie algebra of a link}
\label{subsect:hololink}
Let $K=(K_1,\dots ,K_n)$ be a tame, oriented link in the $3$-sphere,
with complement $X=S^3\setminus \bigcup_{i=1}^n K_i$, and
fundamental group $G=\pi_1(X)$.
By Alexander duality,  $H^1(X)=\Z^n$, generated by classes
$e_1,\dots,e_n$ dual to the meridians of $K$.
As is well-known, the cup-product  map
$\cup_X\colon H^1(X)\wedge H^1(X)\to H^2(X)$ is determined by
the linking numbers, $l_{ij} := \lk (K_i, K_j)$. Consequently, 
the holonomy Lie algebra is determined by the linking
numbers of $K$. More precisely:
\begin{equation}
\label{eq:hlielink}
\HH(X)= \HH(G)= \LL(y_1,\dots ,y_n)\Big/\bigg(\sum_{j=1}^n l_{ij}[y_i,y_j]=0,
\ 1\le i<n\bigg)\, .
\end{equation}

The information coming from linking numbers is conveniently encoded in the 
{\em linking graph} of $K$, denoted by $\G_K$.  This is 
the subgraph of the complete graph on vertices
$\{1,\dots , n\}$, having an edge $\{i,j\}$ whenever $l_{ij} \ne 0$,
weighted by the corresponding linking number $l_{ij}$.

\subsection{Chen groups of $2$-component links}
\label{subsect:chenmura}
In \cite{Mu}, Murasugi computed the Chen groups of a $2$-component link
group $G$. Set $\ell :=\abs{\lk(K_1, K_2)}$.

If $\ell =0$, then $G$ may well be non-$1$-formal. Indeed, in this case 
$\HH(G) =\LL_2$ (cf.~\eqref{eq:hlielink}), and thus, by 
Corollary~\ref{lhilb=cor}, the degree $4$ graded piece of 
$\HH(G)/\HH''(G)$ has rank $3$.  
On the other hand, if $G$ is the group of any $2$-component link, 
then $\theta_4(G)=\phi_4(G)$.  Take $G$ to be the group 
of the Whitehead link.  Then, 
$\gr(G)=\LL(y_1,y_2)/\text{ideal}([y_1, [y_2,[y_1,y_2]]])$, see \cite{Ha}. 
It follows that $\theta_4(G)=2$.  Hence, by Theorem~\ref{thm:hatformal}, 
the fundamental group of the Whitehead link is not $1$-formal.

If $\ell\ne 0$, then $\gr_k(G/G'')=\Z_{\ell}^{k-1}$,  for all $k>1$.
In particular, if $\ell = 1$, then $\gr_{>1}(G/G'')=0$, and
this easily implies $\gr_{>1}(G)=0$, too.

\subsection{The rational Murasugi conjecture}
\label{subsect:qmurasugi}
In view of the above computations, Murasugi made the following
conjecture: If $G$ is the group of an $n$-component link with
all linking numbers equal to $\pm 1$, then, for all $k>1$,
\begin{equation}
\label{eq:muraconj}
\gr_k(G)=\gr_k(\FF_{n-1})\quad\text{and}\quad
\gr_k(G/G'')=\gr_k(\FF_{n-1}/\FF''_{n-1}).
\end{equation}

Murasugi's conjecture was proved, in a more general form,
by Massey and Traldi \cite{MT}. These authors gave a sequence
of necessary and sufficient conditions under which \eqref{eq:muraconj}
holds, including the $\Z$-analogue of condition \eqref{qmb} below.   
Further conditions (for various choices of coefficients) were given 
in \cite{L} and \cite{BP1}.
The following theorem is a rational version of the Murasugi conjecture,
with one additional implication, that takes into account the Lie
algebra structure on the direct sum of the Chen groups.

\begin{theorem}
\label{thm:qmura}
Let $K$ be a link of $n$ components, with complement $X$,
group $G=\pi_1(X)$, holonomy Lie algebra $\HH=\HH(X)$, and
linking graph $\G=\G_K$. The following are equivalent:
\begin{alphenum}
\item \label{qma}
The linking graph $\G$ is connected.

\item \label{qmb}
The rational cup-product map,
$\cup^{\Q}_X\colon H^1(X; \Q)\wedge H^1(X; \Q)\to H^2(X; \Q)$, is onto.

\item \label{qmc}
$\phi_k (G) = \phi_k (\FF_{n-1})$, for all $k > 1$.

\item \label{qmd}
$\theta_k (G) = \theta_k (\FF_{n-1})$, for all $k > 1$.

\item \label{qme}
$\phi_2 (G) = \phi_2 (\FF_{n-1})$.
\end{alphenum}
Moreover, any one of the above conditions implies that:
\begin{itemize}
\item[{\rm (f)}] \label{qmf}
$\gr (G/G'') \otimes \Q = (\HH/\HH'')\otimes \Q$,
as graded Lie algebras.
\end{itemize}
\end{theorem}

\begin{proof} The equivalence \eqref{qma} $\Leftrightarrow$ \eqref{qmb}
is proved  in \cite{L}, the implication \eqref{qma} $\Rightarrow$ \eqref{qmd}
is proved in \cite{MT}, \eqref{qma} $\Rightarrow$ \eqref{qmc} is
proved in \cite{BP1}, while \eqref{qmc} $\Rightarrow$ \eqref{qme} and
\eqref{qmd} $\Rightarrow$ \eqref{qme} are obvious.

To prove \eqref{qme} $\Rightarrow$ \eqref{qmb}, set
$\LL^*_n=\LL(y_1,\dots ,y_n)\otimes \Q$.
Let $\pi \colon \LL_n^2 \surj \gr_2(G) \otimes \Q$ be the
canonical surjection, induced by the Lie bracket of 
$\gr(G)\otimes \Q$. From Sullivan~\cite{Su}, we know 
that $\ker (\pi) =\im (\partial^{\Q}_X)$, where  
$\partial^{\Q}_X$ is the dual of $\cup^{\Q}_X$.
Moreover, $\dim \LL^2_n=\binom{n}{2}$, while 
$\rank (\gr_2(G))=\phi_2(G)=\phi_2(\FF_{n-1})=\binom{n-1}{2}$,
by assumption \eqref{qme}. Thus, $\dim (\ker \pi)=n-1$, which
equals $\dim H^2(X; \Q)$.

We are left with proving \eqref{qma} $\Rightarrow$ (f).
By Lemma 4.1 and Theorems 3.2 and 4.2 from \cite{BP1}, 
condition \eqref{qma} implies
$\gr(G)\otimes \Q =\HH \otimes \Q =\LL_{n-1}\rtimes \LL_1$.
Thus, $\HH'\otimes \Q =\HH^{>1}\otimes \Q =\LL_{n-1}^{>1}=\LL'_{n-1}$, 
and so $\HH''\otimes \Q =\LL''_{n-1}$. Consequently,
\begin{equation}
\label{eq:dimgr}
\dim \left(\gr_k\left( \left( \HH/\HH''\right) \otimes \Q\right)\right) =
\dim \left(\gr_k\left( \left( \LL_{n-1}/\LL_{n-1}''\right) \otimes \Q\right)\right),
\quad\text{for $k>1$}.
\end{equation}
Now recall that the map $\Psi_G^{(2)}\colon \HH/\HH''\to \gr(G/G'')$ is a
surjection. Hence, by the equivalent condition \eqref{qmd} and by
\eqref{eq:dimgr}, the map $\Psi_G^{(2)}\otimes \Q$ is an isomorphism.
\end{proof}

\subsection{The integral Murasugi conjecture}
\label{subsect:zmurasugi}
To state an integral version of the previous theorem, we need to recall
the following definition of Anick \cite{A} (see also Labute \cite{L}).
\begin{definition}
\label{def:pconn}
The linking graph $\G$ is {\em connected modulo $p$}, where $p$  is a
prime number, if there is a spanning subtree of $\G$ whose edges have
linking numbers $l_{ij} \neq 0$ modulo $p$.  The graph $\G$
is {\em strongly connected} if it is connected modulo $p$, for
all primes $p$.
\end{definition}

We also need the following definition from \cite[\S3(2'')]{BP1},
slightly paraphrased for our purposes.
\begin{definition}
\label{def:zgen}
An $n$-component link is called {\em $\Z$-generic} if
its $\Z$-holonomy Lie algebra splits as a semi-direct
product $\LL_{n-1} \rtimes \LL_1$:
\begin{equation}
\label{eq:zgenpres}
\HH=\LL(y_1,\dots ,y_n)/\bigg(
[y_i, y_n] + \sum_{j,k <n} c^i_{j,k} [y_j,y_k]=0, \ 1\le i<n\bigg).
\end{equation}
\end{definition}

We now can improve on Theorem~\ref{thm:qmura},  as follows.
\begin{theorem}
\label{thm:zmura}
Let $K$ be a link of $n$ components, with complement $X$,
group $G=\pi_1(X)$, holonomy Lie algebra $\HH=\HH(X)$, and
linking graph $\G=\G_K$. Then:
\begin{alphenum}
\item \label{zma}
The graph $\G$ is strongly connected if and only if
$\gr(G/G'') = (\LL_{n-1}/\LL''_{n-1}) \oplus \LL_1$
(as graded abelian groups).

\item \label{zmb}
If $\G$ is strongly connected, then the Chen groups
are torsion-free and $\gr(G/G'')=\HH/\HH''$ (as graded Lie
algebras).

\item \label{zmc}
The link $K$ is $\Z$-generic if and only if
the Chen Lie algebra splits as a semi-direct product,
$\gr(G/G'')=(\LL_{n-1}/\LL''_{n-1}) \rtimes \LL_1$.
\end{alphenum}
\end{theorem}

\begin{proof}
\eqref{zma}
As shown by Labute \cite{L}, the graph $\G$ is connected modulo $p$
if and only if the cup-product map
$\cup^{\F_p}_X\colon H^1(X; \F_p)\wedge H^1(X; \F_p)\to H^2(X; \F_p)$ 
is onto.  Thus, $\G$ is strongly connected if and only if
$\cup_X\colon H^1(X; \Z)\wedge H^1(X; \Z)\to H^2(X; \Z)$ is onto.
By Massey-Traldi \cite[Theorem 1]{MT}, this happens if and only if
$\gr_{k} (G/G'')=\gr_{k} (\FF_{n-1}/\FF_{n-1}'')$, for all $k\ge 2$.
Now use Corollary~\ref{hichen=cor}.

\eqref{zmb}
If $\G$ is strongly connected, then $\HH/\HH''$ is torsion-free,
by \cite[3.4--3.6]{MT}. On the other hand, by the implication
\eqref{qma} $\Rightarrow$ (f) from Theorem \ref{thm:qmura},
we know that the surjection
$\Psi_G^{(2)}\colon \HH/\HH''\to \gr(G/G'')$ induces
an isomorphism $\Psi_G^{(2)}\otimes \Q \colon (\HH/\HH'')
\otimes \Q\to \gr(G/G'')\otimes \Q$. It follows that
$\Psi_G^{(2)}$ is an isomorphism, and  $\gr(G/G'')$ is torsion-free.

\eqref{zmc}
Suppose $K$ is $\Z$-generic. The semi-direct product 
decomposition $\HH =\LL_{n-1} \rtimes \LL_1$ from \eqref{eq:zgenpres} 
implies that 
$\partial_X \colon H_2(X;\Z)\to H_1(X;\Z)\wedge H_1(X;\Z)$ is a 
split monomorphism. Therefore, $\G$ is strongly connected, and so  
$\gr(G/G'')=\HH/\HH''$, by Part~\eqref{zmb}.  But now  
$\HH/\HH''= (\LL_{n-1}/\LL''_{n-1}) \rtimes \LL_1$, and this proves 
the forward implication. 

Conversely, suppose
$\gr(G/G'')=(\LL_{n-1}/\LL''_{n-1}) \rtimes \LL_1$.  Then, 
in particular $\gr_2(G)=\gr_2(\FF_{n-1})$. By Theorem 1 from
\cite{MT}, $\G$ must be strongly connected; therefore
$\HH/\HH''= (\LL_{n-1}/\LL''_{n-1}) \rtimes \LL_1$.
Now, 
\begin{equation}
\label{zmura1=eq}
 (\LL_{n-1}/\LL''_{n-1}) \rtimes \LL_1 =\LL(y_1,\dots, y_n)/
\left(\I +\LL''(y_1,\dots, y_{n-1})\right) ,
\end{equation}
where $\I$ is the ideal generated by
$[y_i, y_n] + \sum_{j,k <n} c^i_{j,k} [y_j,y_k]$, for $ 1\le i<n$,
by the very definition of a semi-direct Lie product.
On the other hand,
\begin{equation}
\label{zmura2=eq}
\HH/\HH''=\LL(y_1,\dots, y_n) /
\left( \J +\LL''(y_1,\dots, y_n)\right) ,
\end{equation}
where $\J=$ ideal $(\im \ (\partial_X))$. By comparing 
\eqref{zmura1=eq} and~\eqref{zmura2=eq} in degree $2$, 
we obtain the semi-direct product decomposition \eqref{eq:zgenpres} 
of $\HH$.  This completes the proof.
\end{proof}

\begin{remark}
\label{rem:comp}
The additive part \eqref{zma} of the preceding theorem is basically
contained in the main result of Massey and Traldi \cite{MT}. The other two
multiplicative parts---related to the graded Lie algebra
structure of $\gr(G/G'')$---are new.

Note also that the $\Z$-genericity condition in \eqref{zmc} 
is more restrictive than the strong connectivity condition in 
\eqref{zma} and \eqref{zmb}; see \cite[Example 4.4]{BP1}.
\end{remark}

\section{Chen groups of complex hyperplane arrangements}
\label{sect:hyparr}

Let $\A=\{H_1, \dots, H_n\}$ be an arrangement of $n$ hyperplanes
in $\C^{l}$, with complement $X(\A) =\C^{l}\setminus \bigcup_{i=1}^{n} H_i$.
If $\A$ is central and $\rm{d}\A$ is a decone of $\A$, it is well-known that
$\pi_1(X(\A))=\pi_1(X(\rm{d}\A))\times \Z$. Therefore, there is no loss of
generality in assuming, from now on, that $\A$ is central. Set $X =X(\A)$.

The fundamental group, $G=\pi_1(X)$, admits a finite presentation, with
$n$ generators (corresponding to the meridians), and commutator-relators.

The cohomology ring of the complement of a complex hyperplane arrangement
was computed by Brieskorn \cite{Br}, in answer to a conjecture of Arnol'd. 
A presentation for the ring $H^*(X;\Z)$ in terms of the intersection lattice 
$\mathcal{L}(\A)=\{\bigcap_{H\in {\mathcal B}} H \mid {\mathcal B} \subset \A\}$ 
was given by Orlik and Solomon \cite{OS}. 

An important consequence of Brieskorn's theorem is that every
arrangement complement $X$ is a formal space.  In particular,
the group $G$ is $1$-formal. Thus, our results from Sections
\ref{sect:malcev}--\ref{sect:linalex} may be applied to $G$.

The following combinatorial description of $\HH=\HH(X)=\HH(G)$ 
follows from the Orlik-Solomon theorem. 
Let $[n]=\{1,\dots, n\}$ be the set of points of the underlying matroid 
$\mathcal{M}$ of $\A$. Set $\mathcal{D}:=\{ d \subset [n] \mid
\text{$d$ is a line of $\mathcal{M}$} \}$.
Denote by $J$ the abelian group freely generated by
$\{ (d, j) \mid \text{$d\in \mathcal{D}$  and $j\in d$} \}$,
and define a homomorphism
$\partial \colon J \to \LL^2_Y$ (where $Y=\{ y_1, \dots, y_n \}$) by:
\begin{equation}
\label{eq:pararr}
\partial ((d, j)) =\big[ y_j, \sum_{i\in d} y_i \big] .
\end{equation}
Then the holonomy Lie algebra has presentation
\begin{equation}
\label{eq:12pres}
\HH=\LL_Y/ \text{ideal} \left( \im (\partial \colon J \to
\LL_Y^2)\right) .
\end{equation}
The linearized Alexander invariant, $\B_G=B(\HH)$, is a module over 
the polynomial ring $S=\Z[Y]$, with presentation matrix 
\[
\Delta_G^{(1)} = d_3 + \id\otimes \overline{\partial}_G \colon
S\otimes \left( \Big(\bigwedge\nolimits^3 H_1G\Big) \oplus H_2 G \right)
\longrightarrow \bigwedge\nolimits^2 H_1G. 
\] 
(This matrix first appeared in Theorem 4.6 and Remark 4.7 from \cite{CS3}.)

The next result (which is a direct consequence of Theorem~\ref{thm:hatformal})
gives an affirmative answer to a conjecture formulated in
\cite[\S8]{S}, regarding the combinatorial determination of the Chen
ranks of complex hyperplane arrangements.

\begin{theorem}
\label{thm:arrgts}
Let $\A$ be a complex hyperplane arrangement, with complement $X$,
fundamental group $G$, and holonomy Lie algebra $\HH$, with presentation
given by \eqref{eq:12pres} and \eqref{eq:pararr}. Then:
\begin{equation}
\label{eq:chenarr}
\gr(G/G'')\otimes \Q = (\HH/\HH'') \otimes \Q .
\end{equation}
In particular, the rational Chen Lie algebra of the arrangement
is combinatorially determined (as a graded Lie algebra), by the
level $2$ of the intersection lattice, ${\mathcal L}_2 (\A)$.
Hence, the Chen ranks, $\theta_k(G)=\rank \gr_k(G/G'')$,
are combinatorially determined.
\end{theorem}

\begin{remark}
\label{rem:extors}
The torsion-freeness of $\gr(G/G'')$, also conjectured
in \cite[\S8]{S}, is still open. In a stronger form, this would follow
from
the torsion-freeness of $\HH/\HH''$,  see Proposition~\ref{tors}.
We do have examples of arrangements for which the linearized
Alexander invariant $\B_G$ has torsion, but it is not
clear whether that torsion survives in $\gr(G/G'')$.

The possible presence of torsion in $\gr(G)$ (and consequently, in $\HH$,
see again Proposition~\ref{tors})  was first noticed in \cite[Example
10.7]{S}.
On the other hand, if $\A$ is hypersolvable, then  both $\HH$
and $\gr(G)$ are torsion-free, see \cite[\S7]{JP}.
\end{remark}


\begin{thebibliography}{00}

\bibitem{A} D.~Anick,
{\em Inert sets and the {L}ie algebra associated to a group},
J. of Algebra \textbf{111} (1987), 154--165.

\bibitem{BP1} B.~Berceanu, S.~Papadima,
{\em Cohomologically generic $2$-complexes and $3$-dimensional
Poincar\'{e} complexes}, Math. Ann. \textbf{298} (1994), 457--480.

\bibitem{BP2} \bysame,
{\em Moduli spaces for generic low-dimensional complexes},
J. Pure Appl. Alg. \textbf{95} (1994), 1--25.

\bibitem{Bb} N.~Bourbaki,
{\em Groupes et alg\`{e}bres de Lie}, Chapitres 2--3,
Hermann, Paris, 1972.

\bibitem{Br} E.~Brieskorn,
{\em Sur les groupes de tresses}, in:
S\'{e}minaire Bourbaki 1971/72, Lecture Notes in Math., vol.~317,
Springer-Verlag, Berlin, 1973, pp.~21--44.

\bibitem{Ch51} K.~T.~Chen,
{\em Integration in free groups},
Ann. of Math. \textbf{54} (1951), 147--162.

\bibitem{Ch52} \bysame,
{\em Commutator calculus and link invariants},
Proc. Amer. Math. Soc. \textbf{3} (1952), 44--95.

\bibitem{Ch} \bysame,
{\em Extension of $C^{\infty}$ function algebra by integrals 
and {M}alcev completion of $\pi_1$}, Adv. in Math. \textbf{23} 
(1977), 181--210.

\bibitem{CS1} D.~Cohen, A.~Suciu,
{\em The Chen groups of the pure braid group}, in:
The \v{C}ech Centennial:  A Conference on Homotopy Theory,
Contemp. Math., vol.~181, Amer. Math. Soc., Providence, RI,
1995, pp.~45--64.

\bibitem{CS2} \bysame,
{\em {A}lexander invariants of complex hyperplane arrangements},
Trans. Amer. Math. Soc.  \textbf{351} (1999), 4043--4067.

\bibitem{CS3} \bysame,
{\em Characteristic varieties of arrangements},
Math. Proc. Cambridge Phil. Soc. \textbf{127} (1999), 33--53.

\bibitem{DGMS} P.~Deligne, P.~Griffiths, J.~Morgan, D.~Sullivan,
{\em Real homotopy theory of {K}\"{a}hler manifolds},
Invent. Math. \textbf{29} (1975), 245--274.

\bibitem{FR} M. Falk, R. Randell,
{\em The lower central series of a fiber-type arrangement},
Invent. Math. \textbf{82} (1985), 77--88.

\bibitem{F} R.~H.~Fox,
{\em Free differential calculus. \textup{I}}, 
Ann. of Math. \textbf{57} (1953), 547--560.

\bibitem{Ha} R.~Hain,
{\em Iterated integrals, intersection theory and link groups},
Topology \textbf{24} (1985), 45--66; erratum, 
Topology \textbf{25} (1986), 585--586.

\bibitem{HMR} P.~Hilton, G.~Mislin, J.~Roitberg,
{\em Localization of nilpotent groups and spaces},
North-Holland Math. Studies, vol.~15, North-Holland,
Amsterdam-Oxford; Elsevier, New York, 1975.

\bibitem{JP} M.~Jambu, S.~Papadima,
{\em A generalization of fiber-type arrangements and a new
deformation method}, Topology \textbf{37} (1998), 1135--1164.

\bibitem{Ko83} T.~Kohno,
{\em On the holonomy Lie algebra and the nilpotent completion
of the fundamental group of the complement of hypersurfaces},
Nagoya Math. J. \textbf{92} (1983), 21--37.

\bibitem{Ko85} \bysame,
{\em S\'erie de Poincar\'{e}-Koszul associ\'ee aux groupes
de tresses pures}, Invent. Math. \textbf{82} (1985), 57--75.

\bibitem{Kj} S.~Kojima,
{\em Nilpotent completions and {L}ie rings associated to link groups},
Comment. Math. Helvetici \textbf{58} (1983), 115--134.

\bibitem{L70} J.~P.~Labute,
{\em On the descending central series of groups with a single
defining relation}, J. Algebra \textbf{14} (1970), 16--23.

\bibitem{L} \bysame,
{\em The {L}ie algebra associated to the lower central series
of a link group and Murasugi's conjecture}, Proc. Amer. Math. Soc.
\textbf{109} (1990), 951--956.

\bibitem{La} M.~Lazard,
{\em Sur les groupes nilpotents et les anneaux de {L}ie},
Ann. Sci. \'{E}cole Norm. Sup. \textbf{71} (1954), 101--190.

\bibitem{Md} T.~Maeda,
{\em Lower central series of link groups}, Ph.D. thesis, University
of Toronto, 1983.

\bibitem{MKS} W.~Magnus, A.~Karrass, D.~Solitar,
{\em Combinatorial group theory} (2nd ed.), Dover, New~York, 1976.

\bibitem{MP} M.~Markl, S.~Papadima,
{\em Homotopy {L}ie algebras and fundamental groups via deformation
theory}, Annales Inst. Fourier \textbf{42} (1992), 905--935.

\bibitem{Ms}  W.~S.~Massey,
{\em Completion of link modules}, Duke Math. J.
\textbf{47} (1980), 399--420.

\bibitem{MT}  W.~S.~Massey, L.~Traldi,
{\em On a conjecture of {K}.~{M}urasugi}, Pacific J. Math.
\textbf{124} (1986), 193--213.

\bibitem{MS} D.~Matei, A.~Suciu,
{\em Cohomology rings and nilpotent quotients of real and
complex arrangements}, in: Arrange\-ments--Tokyo 1998, 
Adv. Stud. Pure Math., vol.~27,  Math. Soc. Japan,
Tokyo, 2000, pp.~185--215.

\bibitem{Mo} J.~Morgan,
{\em The algebraic topology of smooth algebraic varieties},
Inst. Hautes \'{E}tudes Sci. Publ. Math. \textbf{48} (1978), 137--204.

\bibitem{Mu}  K.~Murasugi,
{\em On {M}ilnor's invariants for links. {\rm II}.  The {C}hen groups},
Trans. Amer. Math. Soc. \textbf{148} (1970), 41--61.

\bibitem{OS} P.~Orlik, L.~Solomon,
{\em Combinatorics and topology of complements of
hyperplanes}, Invent. Math. \textbf{56} (1980), 167--189.

\bibitem{P} S.~Papadima,
{\em Generalized $\overline{\mu}$--invariants for links and 
hyperplane arrangements},  Proc. London. Math. Soc. \textbf{84} 
(2002), 492--512.

\bibitem{PY}  S.~Papadima, S.~Yuzvinsky,
{\em On rational $K[\pi,1]$ spaces and {K}oszul algebras},
J. Pure Appl. Alg. \textbf{144} (1999), 157--167.

\bibitem{Qu68} D.~Quillen,
{\em On the associated graded ring of a group ring}, J. Algebra
\textbf{10} (1968), 411--418.

\bibitem{Qu} \bysame,
{\em Rational homotopy theory}, Ann. of Math.
\textbf{90} (1969), 205--295.

\bibitem{SS} H.~Schenck, A.~Suciu,
{\em Lower central series and free resolutions of
hyperplane arrangements}, Trans. Amer. Math. Soc.
\textbf{354} (2002), 3409--3433.

\bibitem{S} A.~Suciu,
{\em Fundamental groups of line arrangements: {E}numerative
aspects}, in:  Advances in algebraic geometry motivated by physics, 
Contemp. Math., vol.~276, Amer. Math. Soc, Providence,
RI, 2001, pp. 43--79.

\bibitem{Su} D.~Sullivan,
{\em Infinitesimal computations in topology},
Inst. Hautes \'{E}tudes Sci. Publ. Math.
\textbf{47} (1977), 269--331.

\bibitem{Tr89} L. Traldi,
{\em Linking numbers and {C}hen groups},
Topology Appl. \textbf{31} (1989), 55--71.

\end{thebibliography}
\end{document}